\documentclass{article}
\usepackage{arxiv}
\usepackage{amsfonts,amssymb,amsmath,amsthm}
\usepackage{graphicx,float,multirow}
\usepackage{xcolor}
\usepackage{latexsym,bm}        
\graphicspath{ {./fig/} }

\usepackage{hyperref}
\makeindex

\date{}

\title{Linearly implicit exponential integrators for damped  Hamiltonian PDEs}

\author{
  \href{https://orcid.org/0000-0001-5262-063X}{\includegraphics[scale=0.06]{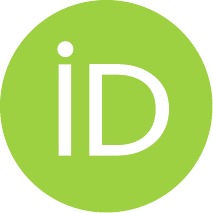}\hspace{1mm}Murat Uzunca} \\
   Department of Mathematics\\
	Sinop University, Sinop-Turkey\\
     \texttt{muzunca@sinop.edu.tr}\\
		    \And
     \href{https://orcid.org/0000-0003-1037-5431}{\includegraphics[scale=0.06]{orcid.pdf}\hspace{1mm}B\"ulent Karas\"ozen} \\
     Department of Mathematics\\
     Middle East Technical University, Ankara-Turkey\\
     \texttt{bulent@metu.edu.tr}
}

\begin{document}

\maketitle

\begin{abstract}
Structure-preserving linearly implicit exponential integrators are constructed for Hamiltonian partial differential equations with linear constant damping. Linearly implicit integrators are derived by polarizing the polynomial terms of the Hamiltonian function and portioning out the nonlinearly of consecutive time steps. They require only a solution of one linear system at each time step. Therefore they are computationally more advantageous than implicit integrators. We also construct an exponential version of the well-known one-step Kahan's method by polarizing the quadratic vector field. These integrators are applied to one-dimensional damped Burger's, Korteweg-de-Vries, and nonlinear Schr{\"o}dinger equations.  Preservation of the dissipation rate of linear and quadratic conformal invariants and the Hamiltonian is  illustrated by numerical experiments.
\end{abstract}

\keywords{
Hamiltonian systems, linear damping, linearly implicit integrator, exponential integrator,  polarization, dissipation preservation}

\section{Introduction}

Many physical systems are commonly affected by external forces or by the dissipative effects of friction. Many well-known partial differential equations (PDEs)  can be expressed with  damping terms:   Burger's equation
 \cite{Bhatt16},  Korteweg-de-Vries (KdV) equation \cite{Bhatt21,Guo19},  nonlinear Schr{\"o}dinger (NLS) equation \cite{Bhatt19,Cai20,Fu16,Moore13},  Klein Gordon equation
\cite{Bhatt16,Moore13},  semi-linear wave equation \cite{Moore13} and  Camassa-Holm equation \cite{Bhatt19} are some known examples.  
In this paper, we consider a Hamiltonian PDE with linear constant damping
\begin{equation}\label{disham} 
u_t  = {\mathcal S } \frac{\delta {\mathcal H}}{\delta u}  -\gamma u ,
\end{equation}
where $u\in \mathbb{R}^d$ is the solution vector for some integer $d\geq 1$, $\mathcal S$ is a constant skew-adjoint operator, $\mathcal H$ is the Hamiltonian, $\delta {\mathcal H}/\delta u$ is the variational derivative and the parameter $\gamma$ stands for the constant damping rate \cite{Bhatt21,Bhatt16,Bhatt17,Moore21}.

Damped PDEs \eqref{disham}  are characterized by the possession of
qualitative properties that decay exponentially along any solution, which are referred to as conformal invariants. 
A conformal invariant ${\mathcal I}(t):\mathbb{R}^d \rightarrow \mathbb{R}$ depending on the solution $u(x,t)$ is defined as \cite{Bhatt16,Moore21,Moore13}
\begin{equation}\label{conf}
\frac{d}{dt} {\mathcal I } = -2\gamma {\mathcal I}, \hbox{ or equivalently }  {\mathcal I} (t) =
e^{-2\gamma t }{\mathcal I} (0).
\end{equation}
This decay in the solution or qualitative
properties of a PDE is often the result of the presence of resistive forces in the system.

The conformal invariant ${\mathcal I}$ describes a quantity of the system such as energy (or Hamiltonian), mass or momentum that decreases with time. This decay in the solution or qualitative properties of a PDE is often due to resistive forces in the system such as friction, damping, dissipation, or viscosity, and hence, are a more realistic model of a physical phenomenon than the conservative systems.  

It is important to maintain as many properties of the physical system as possible when modeling physical phenomena of some useful discretization techniques. 
Numerical methods, especially the energy-preserving methods, for conservative and dissipative systems, have attracted a significant amount of attention in recent years. 
Numerical methods that preserve the conformal symplectic structure of conformal Hamiltonian systems are known as conformal symplectic methods. 
They were first constructed for ordinary differential equations (ODEs) using splitting techniques \cite{mclachlan02} by solving the linear dissipative part exactly and the nonlinear conservative part with a symplectic method, and then composing the flow
maps. 
Various integrators are constructed using splitting techniques preserving the conformal multi-symplectic structure of damped PDEs; the KdV equation \cite{Guo19}, 
the NLS equation \cite{Fu16}, semi-linear wave equation \cite{Moore09}. 
Other conformal structure-preserving integrators are conformal multi-symplectic Euler-Preissman scheme \cite{Moore17} and discrete gradient method \cite{Moore13},  St{\"o}rmer-Verlet and conformal implicit midpoint methods \cite{Bhatt16}, exponential Rung-Kutta methods \cite{Bhatt19}, projected exponential Runge-Kutta methods \cite{Bhatt21}.

In this paper,  we construct the linearly implicit exponential integrators for damped Hamiltonian PDEs \eqref{disham} by combining the linearly implicit methods using polarized energy \cite{Dahlby11,Li22}  with the exponential methods using discrete gradient \cite{Moore21}. Implicit exponential integrators are constructed for 
damped PDEs \eqref{disham} using discrete gradients in \cite{Moore21} such as the exponential average vector field method and exponential implicit midpoint method.
They can be considered as extension of the energy-preserving discrete gradient methods for Hamiltonian PDEs through the development of exponential integrators. 
Some numerical
methods preserve the dissipation properties by simply guaranteeing that the energy or conformal invariant is
decreasing with every time step, even though it may be numerically overdamped
or underdamped.
The exponential integrators in \cite{Bhatt21,Bhatt16,Moore21} preserve the correct rate of dissipation, such that the energy or the conformal invariant  is not overdamped or
underdamped. 
Due to their implicit nature, a system of nonlinear equations have to be solved iteratively at each time step by Newton's method or by fixed point iteration.
On the other hand, linearly implicit integrators require only a single iteration in the solution of a nonlinear system of equations, which makes the linearly implicit integrators computationally advantageous over the implicit exponential integrators such as the average vector field (AVF) and conformal midpoint methods. 
The linearly implicit methods are constructed using polarized energy to portion out the nonlinear terms in the Hamiltonian function over consecutive time steps. In this way, a quadratic polarized energy is constructed and then the polarized discrete gradient method is performed. Linearly implicit energy-preserving methods have been applied to  Hamiltonian PDEs  \cite{Cai22,Eidnes20,Eidnes21} and gradient systems 
\cite{Karasozen23ac} with polynomial nonlinear terms.  
In this study, we derive two-step linearly implicit exponential methods for damped Burger's, KdV, and NLS equations where the Hamiltonian functions contain quadratic, cubic, and quartic terms. 
Linearly implicit methods are symmetric, preserve the polarized energy, and have favorable properties like linear error growth and long-time near-conservation of first integrals.  
Similarly, linearly implicit exponential integrators also preserve the correct dissipation rate of the Hamiltonians and conformal invariants.
A well-known one-step linearly implicit integrator for linear-quadratic systems is  Kahan's method \cite{Kahan93,Kahan97} which is constructed by polarizing the quadratic vector fields. We construct  linearly implicit exponential two-step Kahan's method for damped Burger's and KdV equations \cite{Celledoni13}. 
Preservation of the dissipation of linear and quadratic conformal invariants and the energy are illustrated for damped Burger's, KdV, and NLS equations through numerical experiments

The rest of this paper is organized as follows. In Section~\ref{sec2}, linearly implicit exponential integrators are introduced. Numerical results for damped Burger's, KdV  and NLS equations are presented in Section~\ref{sec:num}. The paper ends with some conclusions in Section~\ref{sec:conc}.

\section{Linearly implicit exponential integrators}
\label{sec2}

Semidiscretization of  \eqref{disham} in space with finite differences, yields the  dissipative Hamiltonian ODE  
\begin{equation} \label{dishamode}
\dot{\bm u} = S \nabla H (\bm u) - \gamma{\bm u},
\end{equation}
where ${\bm u}(t): [0,T]\longrightarrow \mathbb{R}^M $ is the unknown solution vector, $M$ is the spatial degrees of freedom, and $S$ is  $M\times M$ constant skew-symmetric matrix.
It is desirable that a numerical method to the semidiscrete system \eqref{dishamode}, preserve the conformal invariant \eqref{conf} numerically
\begin{equation}
I^{n+1} = e^{-2\gamma\Delta t} I^{n},\quad n=0,1,\ldots
\end{equation}
where $I^n\approx {\mathcal I}(t_n)$ is a discrete approximation of the conformal invariant \eqref{conf}  at time $t_n$, and $\Delta t$ is the time step size.
Several methods have been developed for preserving dissipation of the conformal invariants  \cite{Bhatt21,Bhatt16,Bhatt17,Guo19,Moore13}. In this paper, we construct linearly implicit exponential discrete gradient methods for \eqref{dishamode} following \cite{Moore21,Eidnes21,Dahlby11}. In the first part, we briefly describe discrete gradient methods and linearly implicit 
integrators using quadratic polarization of the Hamiltonians. In the second part, we derive the linearly implicit  exponential integrators for dissipatively perturbed Hamiltonian systems.

\subsection{Discrete gradient methods and linearly implicit integrators}

For conservative Hamiltonian systems without damping $(\gamma = 0$)
\begin{equation}\label{sham}
\dot{\bm u} = \bm{f} ( {\bm u}  ) = S \nabla H (\bm u) ,
\end{equation}
the discrete gradient method is given by
\begin{equation}\label{dgrad}
\frac{{\bm u}^{n+1}-{\bm u}^n}{\Delta t}=S\overline{\nabla}H({\bm u}^n,{\bm u}^{n+1}),
\end{equation}  
where ${\bm u}^{n}={\bm u}(t_n)$ is the solution vector approximating $u(t_n)$.

The discrete gradient $\overline{\nabla}H : \mathbb{R}^M{\times} \mathbb{R}^M\rightarrow \mathbb{R}^M$ is a vector such that for any vectors ${\bm u},{\bm v} \in \mathbb{R}^M$, it holds
\begin{itemize}
\item $\displaystyle H({\bm v})-H({\bm u})= ({\bm v}-{\bm u})^T\overline{\nabla}H({\bm u},{\bm v})$,
\item $\displaystyle \overline{\nabla}H({\bm u},{\bm u}) =   \nabla H({\bm u})$.
\end{itemize}
The discrete gradient method \eqref{dgrad} preserves the energy of the conservative Hamiltonian system  \eqref{sham}  at any time step. 
A particular choice of $\overline{\nabla}H$ such as mean value approximation
\begin{equation*} 
\overline{\nabla}H({\bm u}^n,{\bm u}^{n+1}) =  
\int_0^1 \nabla H \left ( \xi {\bm u}^{n+1} + (1-\xi) {\bm u}^{n} ) d\xi \right ) ,
\end{equation*}
leads to the AVF method for Hamiltonian systems \cite{Celledoni12,Cohen11,Quispel08}.
These methods are implicit, which means a nonlinear system of equations has to be solved iteratively at each time step. 
Linearly implicit integrators require the solution of one linear system
of equations at each time step.
Linearly implicit methods are constructed by portioning out the nonlinearity over consecutive time steps by devising 
 a quadratic polarized energy $\widetilde{H}({\bm v},{\bm w }):\mathbb{R}^M{\times} \mathbb{R}^M\rightarrow \mathbb{R}^M$ satisfying consistency and invariance properties
\begin{equation}\label{psym}
\begin{aligned}
\widetilde{H}({\bm u}, \bm{u}) &= H({\bm u}), \\
\widetilde{H}({\bm u}, \bm{v}) &=\widetilde{H}({\bm v}, \bm{u}), 
\end{aligned}
\end{equation}
and then performing the polarized discrete gradient (PDG) method.

A PDG for $\widetilde{H}$ is a function $\overline{\nabla}\widetilde{H} : \mathbb{R}^M \times \mathbb{R}^M\times \mathbb{R}^M \rightarrow \mathbb{R}^M$ satisfying \cite{Cai22,Dahlby11,Eidnes20,Eidnes21,Li22,Karasozen23ac}
\begin{equation}\label{dg1}
\begin{aligned}
\widetilde{H}({\bm v},{\bm w}) - \widetilde{H}({\bm u},{\bm v}) &= 
\frac{1}{2}({\bm w}-{\bm u})^T \overline{\nabla}\widetilde{H}({\bm u},{\bm v},{\bm w} ), \\
\overline{\nabla}\widetilde{H}({\bm u},{\bm u}, {\bm u })&= \nabla H({\bm u }).
\end{aligned}
\end{equation}
The corresponding polarized two-step discrete gradient scheme is given by  \cite{Dahlby11,Li22}
\begin{equation}\label{twstepli}
\frac{{\bm u}^{n+2}-{\bm u}^{n}}{2 \Delta t} = S \overline{\nabla}\widetilde{H}({\bm u}^n,{\bm u}^{n+1},{\bm u}^{n+2}),
\end{equation}
which preserves the polarized Hamiltonian $\widetilde{H}$ in the sense $\widetilde{H}({\bm u}^{n},{\bm u}^{n+1}) =  \widetilde{H}({\bm u}^{0},{\bm u}^{1})$ for all $n\ge 0$ \cite{Eidnes21}.
Quadratic polarization of several for polynomial functions are given below  \cite{Dahlby11,Li22}:
\begin{itemize} 
    \item  $\displaystyle H(\bm u)={\bm u}^2  \rightarrow \;\widetilde{H}({\bm v},{\bm w})=\theta\frac{{\bm v}^2+{\bm w}^2}{2}+(1-\theta){\bm v}{\bm w}$, \quad $\theta\in [0,1]$,
    \item   $\displaystyle H({\bm u})={\bm u}^3 \;   \rightarrow  \; \widetilde{H}({\bm v},{\bm w})={\bm v}{\bm w}\frac{{\bm v}+{\bm w}}{2}$,
      \item   $\displaystyle H({\bm u})={\bm u}^4 \; \rightarrow  \; \widetilde{H}({\bm v},{\bm w})={\bm v}^2{\bm w}^2$.
\end{itemize}

For linear-quadratic systems such as Burger's equation and KdV equation with cubic Hamiltonian functions, there exist a well known linearly implicit method, namely Kahan's method \cite{Celledoni13,Kahan93,Kahan97}. 
When restricted to quadratic vector fields, Kahan's method coincides with the Runge-Kutta method \cite{Celledoni13}
\begin{equation} \label{kahan}
\frac{{\bm u}^{n+1} - {\bm u}^n } {\Delta t} = -\frac{1}{2}\bm{f}({\bm u}^{n+1})   +  2\bm{f}\left( \frac{{\bm u}^{n+1} + {\bm u}^n }{2}  \right)   -\frac{1}{2}\bm{f}({\bm u}^{n} ).
\end{equation}

The two-step linearly implicit PDG scheme \cite{Eidnes21} 
\begin{equation}\label{eq:kahan2s}
\frac{{\bm u}^{n+2}-  {\bm u}^n}{2\Delta t} = \frac{1}{4}S H''({\bm u}^{n+1})({\bm u}^n+{\bm u}^{n+2}),
\end{equation}
is equivalent to Kahan's method \eqref{kahan} over two consecutive steps, when applied to ODEs with homogeneous cubic $H$
\begin{equation}\label{Two_step_kahan}
\frac{{\bm u}^{n+2}-{\bm u}^n}{2\Delta t}=-\frac{1}{4}\bm{f} ({\bm u}^{n})-\frac{1}{2}\bm{f} ({\bm u}^{n+1})-\frac{1}{4} \bm{f}({\bm u}^{n+2})+ \bm{f}(\frac{{\bm u}^{n}+ 
{\bm u}^{n+1}}{2})+ \bm{f}(\frac{{\bm u}^{n+1}+ {\bm u}^{n+2}}{2}).
\end{equation}
The two-step PDG scheme \eqref{Two_step_kahan} preserves the polarized invariant \cite{Eidnes21}
\begin{equation}\label{polhamk} 
\tilde{H}({\bm u}^n,{\bm u}^{n+1})= \frac{1}{6}({\bm u}^n)^T H''(\frac{{\bm u}^n+{\bm u}^{n+1}}{2}){\bm u}^{n+1},
\end{equation} if  the one-step Kahan's method \eqref{kahan} is used to calculate $\bm u^{1}$ from $\bm u^0$. 
The one-step \eqref{kahan}   and two-step \eqref{Two_step_kahan} Kahan's methods are time-symmetric (reversible) and therefore second-order accurate in time.

\subsection{Linearly implicit exponential integrators}

For dissipative Hamiltonian systems \eqref{disham},  in \cite{Moore21}, the following exponential discrete gradient method  was introduced 
\begin{equation*}
\frac{e^{X_{1}}{\bm u}^{n+1} - e^{X_0}{\bm u}^n  }{\Delta t }= S \overline{\nabla}H( e^{X_0}{\bm u}^n, e^{X_1}{\bm u}^{n+1}),
\end{equation*}
which preserves some particular choices of the conformal invariant \eqref{conf}, where 
$
X_{\alpha} =\int_{t_{n+1/2}}^{t_\alpha}\frac{\gamma}{2}dt. 
$
The mean value averaged  exponential discrete gradient method  leads to the exponential AVF method \cite{Moore21}
\begin{equation} \label{confavf}
\frac{e^{X_{1}}{\bm u}^{n+1} - e^{X_0}{\bm u}^n  }{\Delta t } = S 
\int_0^1 \nabla H \left ( \xi e^{X_{1}}{\bm u}^{n+1} + (1-\xi) e^{X_{0}}{\bm u}^{n}  \right) d\xi,  
\end{equation}
with the  Hamiltonian preserved along the transformed solution \cite{Moore21}
$$
H(e^{X_1}{ \bm u} ^{n+1}) = H(e^{X_0}{\bm u}^{n}), \quad n =0,1,\ldots
$$

The conformal implicit midpoint method \cite{Bhatt16,Moore21}
\begin{equation} \label{confmid}
\frac{e^{X_{1}}{\bm u}^{n+1} - e^{X_0}{\bm u}^n  }{\Delta t } = S \nabla H \left ( \frac{e^{X_{1}}{\bm u}^{n+1}+ e^{X_{0}}{\bm u}^{n}}{2} \right ),
\end{equation}
is equivalent to the mean value averaged  exponential discrete gradient method \eqref{confavf} for linearly damped Hamiltonian systems with cubic Hamiltonian functions such as the damped Burger's equation  and damped KdV equation.  It is symmetric and hence second-order 
\cite{Moore21}.

The exponential version of the one-step exponential Kahan's method \eqref{kahan} can be given as 
\begin{equation}\label{expkahan}
\frac{e^{X_1}{\bm u }^{n+1} - e^{X_0}{\bm u }^{n} }{\Delta t} = -\frac{1}{2}\bm{f}(e^{X_0}{\bm u }^{n} ) +  2\bm{f}\left( \frac{e^{X_1}{\bm u }^{n+1}  + e^{X_0}{\bm u }^{n} }{2}  \right)   -\frac{1}{2}\bm{f}(e^{X_1}{\bm u }^{n+1} ),
\end{equation} 
whose adjoint equation
\begin{equation*}
\frac{e^{X_0}{\bm u }^{n} - e^{X_1}{\bm u }^{n+1} }{-\Delta t} =  -\frac{1}{2}\bm{f}(e^{X_1}{\bm u }^{n+1}) +  2\bm{f}\left( \frac{e^{X_0}{\bm u }^{n} + e^{X_1}{\bm u }^{n+1}  }{2}  \right) - \frac{1}{2}\bm{f}(e^{X_0}{\bm u }^{n}) .
\end{equation*}
is the same as \eqref{expkahan}.  Hence it is symmetric and  second-order.

Similarly, the two-step exponential Kahan's method has the form
\begin{equation}\label{Two_step_expkahan}
\begin{aligned}
\frac{e^{X_2}{\bm u}^{n+2}-e^{X_0}{\bm u}^n}{2\Delta t}   =  & -\frac{1}{4} \bm{f} (e^{X_0}{\bm u}^{n})-\frac{1}{2} \bm{f} (e^{X_1}{\bm u}^{n+1})-
\frac{1}{4} \bm{f} (e^{X_2}{\bm u}^{n+2})\\ 
& +  \bm{f} \left (\frac{e^{X_0}{\bm u}^{n}+ e^{X_1}{\bm u}^{n+1}}{2}\right )+ \bm{f} \left (\frac{e^{X_1}{\bm u}^{n+1}+ {e^{X_2}\bm u}^{n+2}}{2}\right ),
\end{aligned}
\end{equation}

It is symmetric under $(e^{X_0}{\bm u}^n,e^{X_1}{\bm u}^{n+1},e^{X_2}{\bm u}^{n+2}) \rightarrow  (e^{X_2}{\bm u}^{n+2},e^{X_1}{\bm u}^{n+1},e^{X_0}{\bm u}^{n})$ and $\Delta t  \rightarrow  -\Delta t$ as the right hand side  is coming from a symmetric multilinear form. Therefore the scheme \eqref{Two_step_expkahan} is second-order accurate.

The two-step linearly implicit exponential   discrete gradient methods for  damped Hamiltonian systems have the general form 
\begin{equation}\label{polscheme}
\frac{e^{X_2}{\bm u}^{n+2}-e^{X_0}{\bm u}^{n}}{2 \Delta t} = S \overline{\nabla}\widetilde{H}(e^{X_0}{\bm u}^n,e^{X_1}{\bm u}^{n+1},e^{X_2}{\bm u}^{n+2}),
\end{equation}
where 
$
X_{\alpha} =\int_{t_{n+1}}^{t_\alpha}\frac{\gamma}{2}dt. 
$

Like the exponential Kahan's method \eqref{expkahan}, the exponential polarized two-step discrete gradient method \eqref{polscheme} is also second-order, since the adjoint equation 
\begin{equation*}
\frac{e^{X_0}{\bm u}^{n}-e^{X_2}{\bm u}^{n+2}}{-2 \Delta t} = S \overline{\nabla}\widetilde{H}(e^{X_2}{\bm u}^{n+2},e^{X_1}{\bm u}^{n+1},e^{X_0}{\bm u}^{n}).
\end{equation*}
is the same as the \eqref{polscheme}, since the right-hand side  of \eqref{polscheme} is symmetric.

Both \eqref{Two_step_expkahan} and \eqref{polscheme} preserve the polarized Hamiltonian   \eqref{polhamk} along the transformed solution
$$
\widetilde{H}(e^{X_0}{ \bm u} ^{n},  e^{X_1}{ \bm u} ^{n+1}) = \widetilde{H}(e^{X_0}{ \bm u} ^{0},  e^{X_1}{ \bm u} ^{1}) , \quad  n =0,1,\ldots
$$
Linearly implicit two-step exponential discrete gradient integrators  using the polarization  of quadratic, cubic, and quartic terms in Hamiltonian function for damped Burger's,  damped KdV, and  damped NLS equations are given in Section~\ref{sec:num}.

For the dissipative Hamiltonian systems, the Hamiltonian and the conformal invariants dissipate like
$
 H(\bm u (t) ) = e^{-\gamma t} H ({\bm u}^0). \qquad I (\bm u (t) ) = e^{-\gamma t} I ({\bm u}^0).
$
When they are initially positive, exponential integrators guarantee that the Hamiltonian, and conformal invariants are decreasing at each time step, i.e.
$$
H(t_{n+1}) < H(t_n), \qquad  I(t_{n+1}) < I(t_n), 
$$
for small values of the damping term $\gamma$,  and satisfy energy and conformal invariant balance equations \cite{Bhatt16,Moore21}
\begin{equation*} 
H(e^{X_1}{\bm u}^{n+1})  - H(^{X_0}{\bm u}^{n})=0, \quad I(2e^{X_1}{\bm u}^{n+1})  - I(2e^{X_0}{\bm u}^{n})=0.
\end{equation*}

Preservation of the Hamiltonian (or energy) dissipation is measured using with residual \cite{Moore21} 
\begin{equation} \label{hambal}
R_H = \ln \left( \frac{H(\bm{u}^{n+1})}{H(\bm{u}^{n})}\right) + \Delta t \gamma.
\end{equation}
Similarly, the dissipation of the conformal invariants is measured with the residual \cite{Bhatt16,Moore21}
\begin{equation} \label{mombal}
R_I = \ln \left( \frac{ I (\bm{u}^{n+1})}{ I (\bm{u}^{n})}\right) + 2\Delta t \gamma.
\end{equation}
Both the residuals \eqref{hambal} and \eqref{mombal} is used to check whether the decrease of Hamiltonian and conformal invariants are over/underdamped or not.

\section{Numerical experiments}
\label{sec:num}

In this Section, we report on numerical experiments for three Hamiltonian PDEs with constant linear damping under periodic boundary conditions; Burger's equation, KdV equation, and NLS equation. They are discretized in space with finite differences and the resulting dissipative Hamiltonian ODEs are solved with the implicit and linearly implicit exponential integrators: the conformal implicit midpoint method (CIMP) \eqref{confmid}, the exponential AVF method (EAVF) \eqref{confavf}, the two-step exponential Kahan's method (EK) \eqref{Two_step_expkahan}, and the two-step linearly implicit exponential integrator (LIE) \eqref{polscheme}.

Space discretization on the interval $x \in [-L, L]$ is performed
by introducing a uniform spatial grid of the nodes $\{x_1, x_2, \ldots, x_M\}$ with the grid size $\Delta x$ such that $x_1 =-L, x_M = L$, and $M$ is even. Then, for any $t$, we approximate $u(x_k,t)$ by $u_k(t)$, $k = 1, 2, \ldots, M$, with periodic boundary conditions $u_{k+M}(t) =u_k(t)$, and define the solution vector
${\bm u}(t) = (u_1(t),\cdots,u_M(t))^T$.
The  matrices $D_1\in\mathbb{R}^{M\times M}$ and $D_2\in\mathbb{R}^{M\times M}$ correspond to the centered finite difference discretization of the first and second-order derivative operators $\partial_x$ and $\partial_{xx}$, respectively, under periodic boundary conditions  
\begin{equation*} 
 D_1=\frac{1}{2\Delta x}
	\begin{pmatrix}
	0 & 1 & & & -1\\
	-1 & 0 & 1 & & \\
	  & \ddots &\ddots &\ddots & \\
	 &  & -1 & 0 & 1\\
	1 &  & & -1 & 0\\
	\end{pmatrix},
\quad 
D_2 = \frac{1}{\Delta x ^2}
	\begin{pmatrix}
-2 & 1 & & & 1\\
	1 & -2 & 1 & & \\
	  & \ddots &\ddots &\ddots & \\
	 &  & 1 & -2 & 1\\
	1 &  & & 1 & -2\\
	\end{pmatrix}.
\end{equation*}

For time discretization, we divide the time interval $[0,T]$ into $N$ uniform elements $
0=t_0<t_1 < \cdots < t_{N}=T$, $\Delta t = T/N$, and we denote by $\bm{u}^n= \bm{u}(t_n)$ the full discrete approximation vector at time $t_n$, $n=0,\ldots , N$.

The computations are carried out via MATLAB 7.0 with Intel(R) Core(TM) i5-7500 CPU \@ 3.40 GHz.

\subsection{Damped Burger's equation}

The damped or modified Burger's equation \cite{Bhatt16}
\begin{equation*} 
u_t =  -u u_x- 2\gamma u,
\end{equation*}
can be written as a damped Hamiltonian PDE \eqref{disham}
\begin{equation} \label{burgers}
u_t = -\frac{\partial }{\partial x } \frac{\delta {\mathcal H}}{\delta u }  - 2\gamma u, \qquad {\mathcal H}(u)  = \frac{1}{3}u^3.
\end{equation}
The semidiscrete system reads as
$$
\dot{{\bm u}} = -\frac{1}{2}D_1  ({\bm u} \odot  {\bm u} )-2 \gamma {\bm u},
$$
where $\odot$ denotes the elementwise vector or matrix multiplication.
The Hamiltonian of the full discrete system  is given as 
$$
H({\bm u}^n)  = \frac{1}{3}\int \left (  ( {\bm u}^n )^3 \right ) dx. 
$$

Applying the exponential integrators presented in Section~\ref{sec2}, we obtain the following schemes:
\begin{itemize}
\item CIMP: 
\begin{equation*}
\frac{e^{X_1}{\bm u}^{n+1} - e^{X_0}{\bm u}^n  }{\Delta t }  =   -D_1 \frac{e^{X_0}{\bm u}^{n}\odot e^{X_0}{\bm u}^{n}, 
+    e^{X_1}{\bm u}^{n+1}\odot e^{X_1}{\bm u}^{n+1} }{4 }
\end{equation*}
where
$
 \quad X_{1} = \frac{\gamma\Delta t}{2}, \quad X_{0}= -\frac{\gamma\Delta t}{2}.
$

\item EK: 
\begin{equation*}
\frac{e^{X_2}{\bm u}^{n+2} - e^{X_0}{\bm u}^n  }{2\Delta t } = -D_1\frac{e^{X_1}{\bm u}^{n+1}\odot  e^{X_0}{\bm u}^n  + e^{X_2} {\bm u}^{n+2}\odot  e^{X_1}{\bm u}^{n+1} }{4 },
\end{equation*}
where
$
X_{2} = \frac{\gamma\Delta t}{2},
 \quad X_{1} = 0, \quad X_{0}= -\frac{\gamma\Delta t}{2}.
$

\end{itemize}

Damped Burger's equation \eqref{burgers} possess a  linear conformal invariant,  mass $C(u) = \int udx$ \cite{Bhatt16}.

As the numerical test problem, we consider damped Burger's equation \eqref{burgers} \cite{Bhatt16} on $x\in [-\pi,\pi]$  with the damping factor $\gamma = 0.25$.  Initial condition is taken as the Gaussian distribution with density $1$ and mean zero, i.e., $u(x,0) = e^{-x^2/2}/\sqrt{2\pi}$. We set the spatial and temporal mesh sizes $\Delta x=\pi/40$ and $\Delta t=0.009$, respectively, and the target time $T=50$. We show the numerical solution and decrease of the conformal invariant and of the Hamiltonians only  for the EK. 
The numerical results obtained with CIMP are very close to those obtained with EK, therefore thy are not shown.

\begin{figure}[th]
\centering
\includegraphics[width=.55\columnwidth]{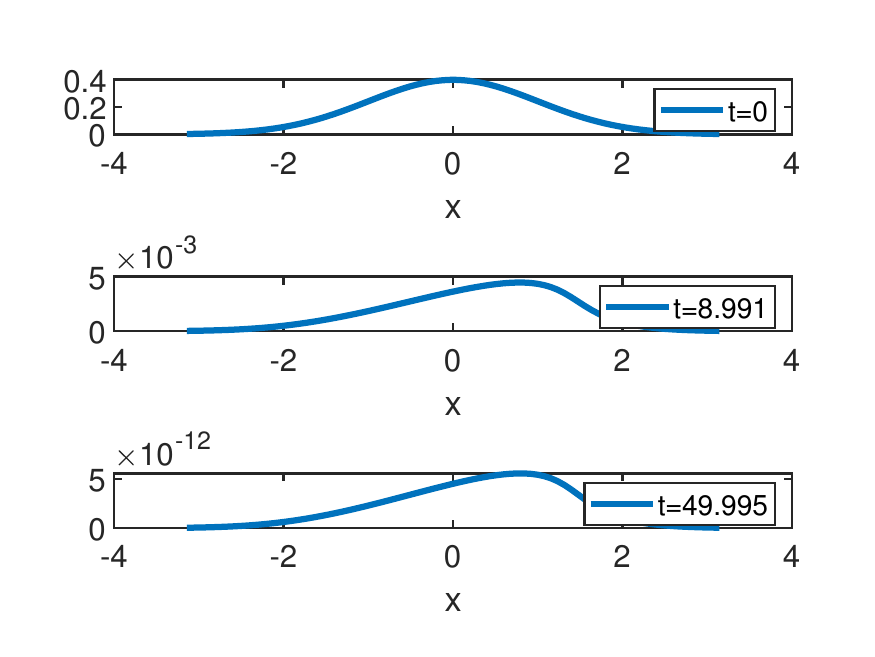}
\caption{Solutions of the damped Burger's equation with EK.\label{burgerssol}}
\end{figure}

\begin{figure}[th]
\centering
\includegraphics[width=.60\columnwidth]{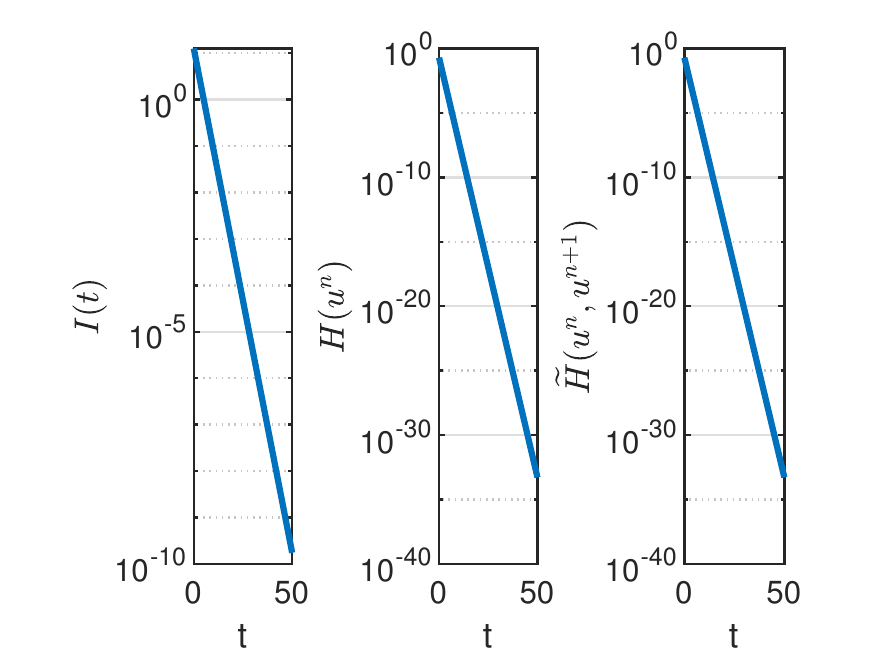}
\caption{Damped Burger's equation (EK):  mass  (left), exact Hamiltonian (middle) and modified Hamiltonian (right).\label{hamburgers}}
\end{figure}

The solution profile becomes steeper and decreases gradually with time in Figure \ref{burgerssol} as in \cite{Bhatt16}. The linear conformal invariant, mass, Hamiltonian, and modified Hamiltonian decrease as time progress in Figure~\ref{hamburgers}. 

\begin{figure}[th]
\centering
\includegraphics[width=.9\textwidth]{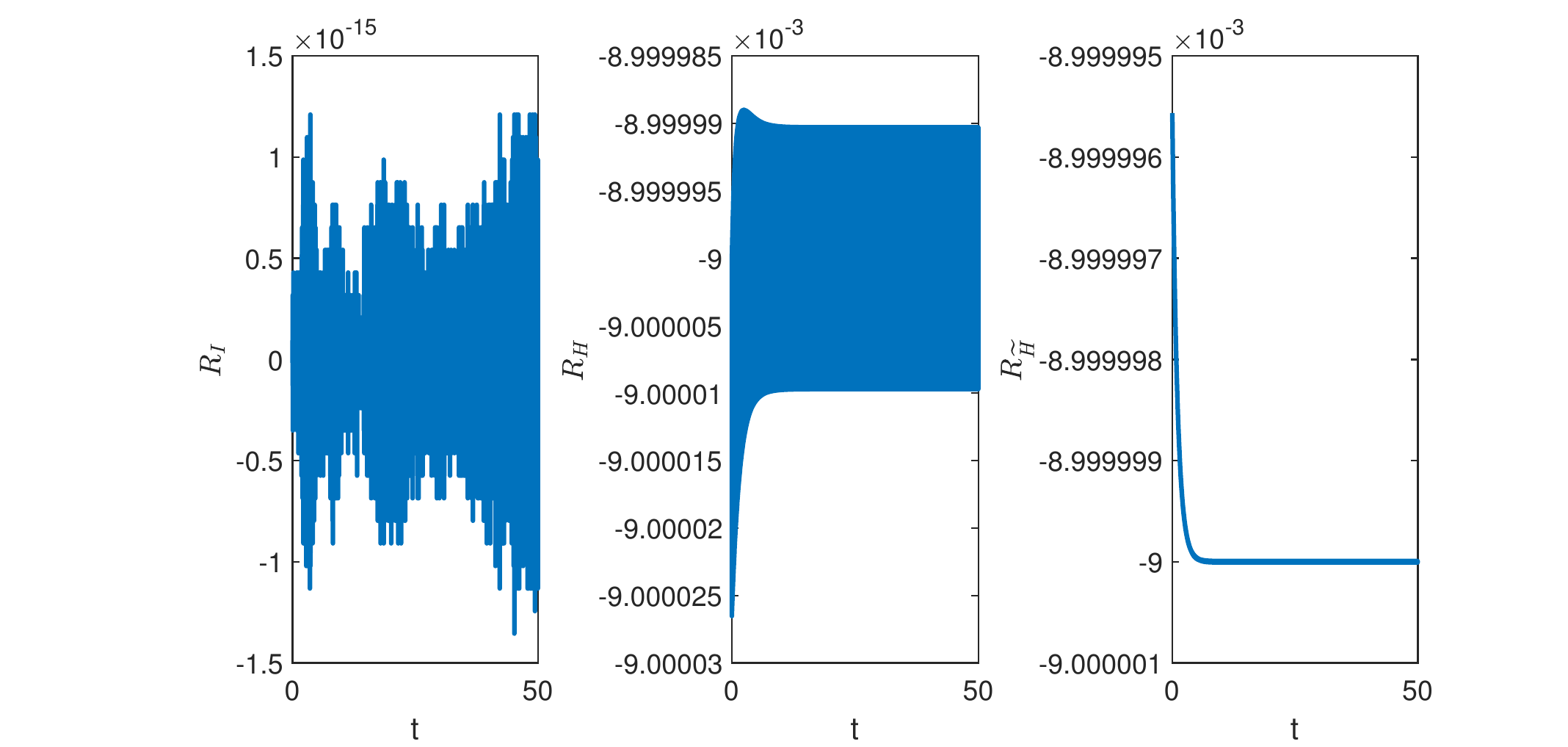}
\caption{Damped Burger's equation (EK):  residuals of mass (left),   Hamiltonian  (middle) and modified Hamiltonian (right).\label{balanceK}}
\end{figure}

\begin{figure}[th]
\centering
\includegraphics[width=.9\textwidth]{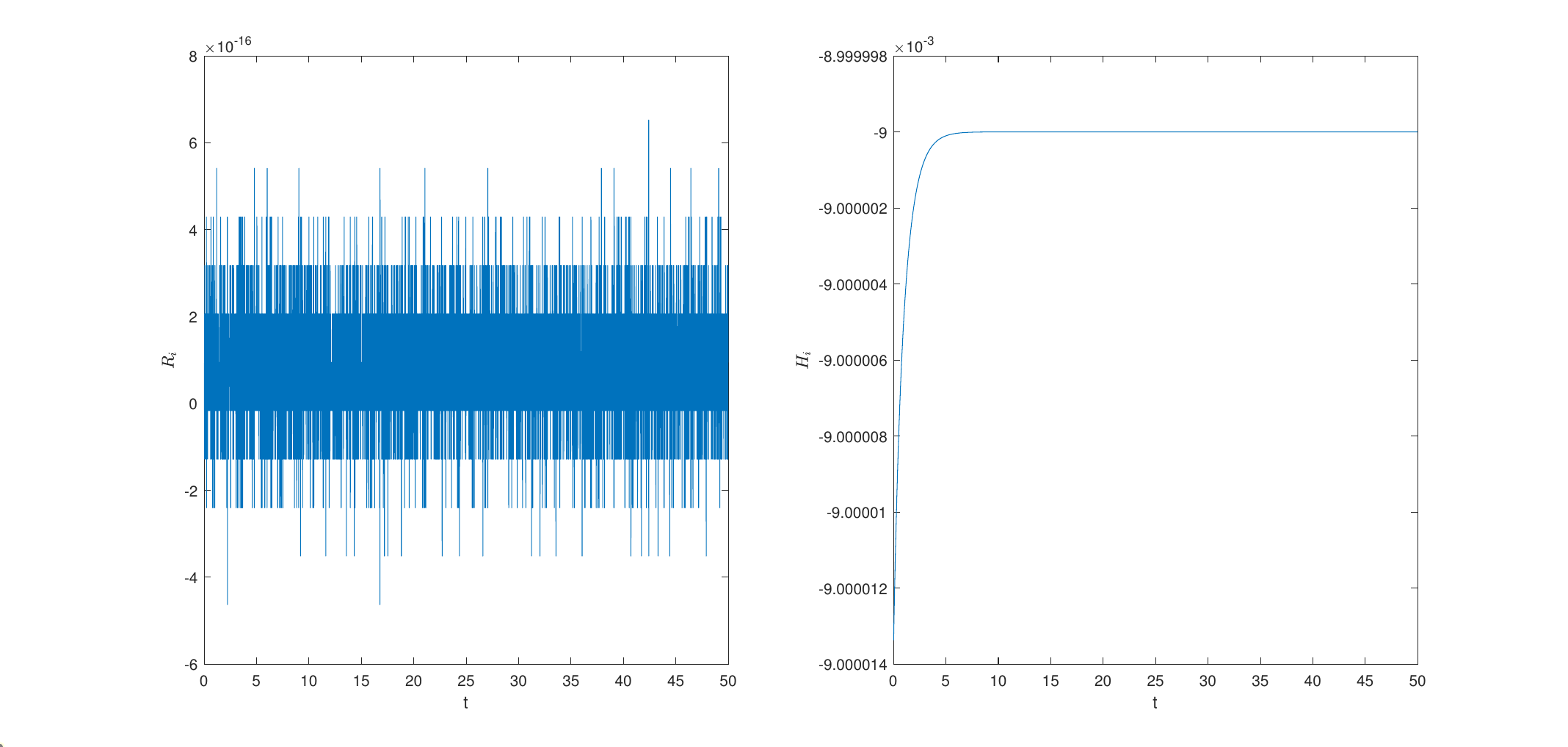}
\caption{Damped Burger's equation (CIMP):  residuals of mass (left),   Hamiltonian  (right).\label{balancecimp}}
\end{figure}

The error in the residual \eqref{mombal} of mass is preserved up to machine precision  in Figure \ref{balanceK}.  The error in the residual  of the Hamiltonians \eqref{hambal} is much larger than of the mass,  but does not show any drift as time progresses, which indicates that they are over-or under damped. 
The mass is also preserved with high accuracy with CIMP as shown in Figure \ref{balancecimp}.

\begin{figure}[th]
\centering
\includegraphics[width=.6\columnwidth]{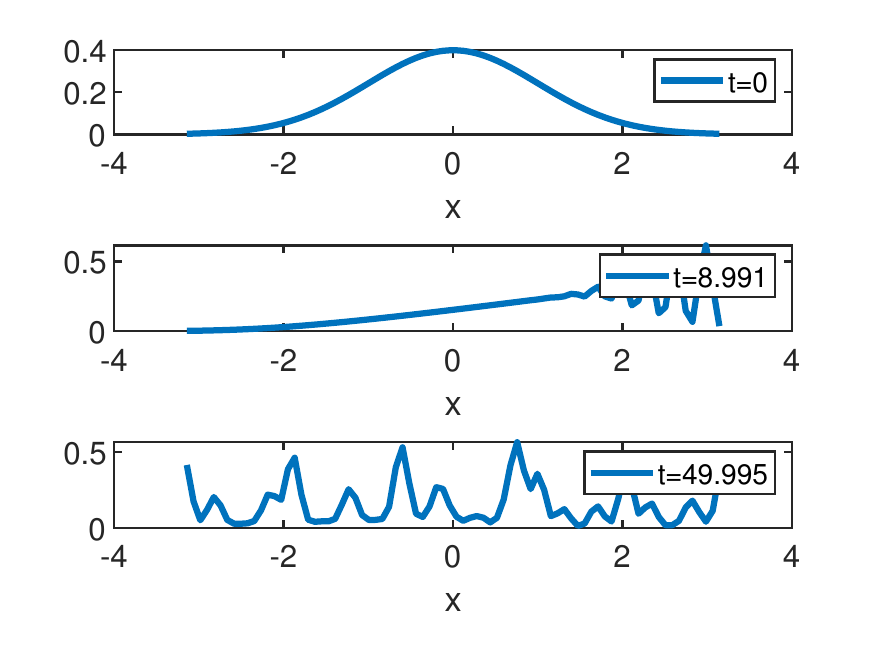}
\caption{Solutions of the damped Burger's equation with Kahan's method \eqref{Two_step_expkahan}.\label{burgerssolkahan}}
\end{figure}

\begin{figure}[th]
\centering
\includegraphics[width=.9\textwidth]{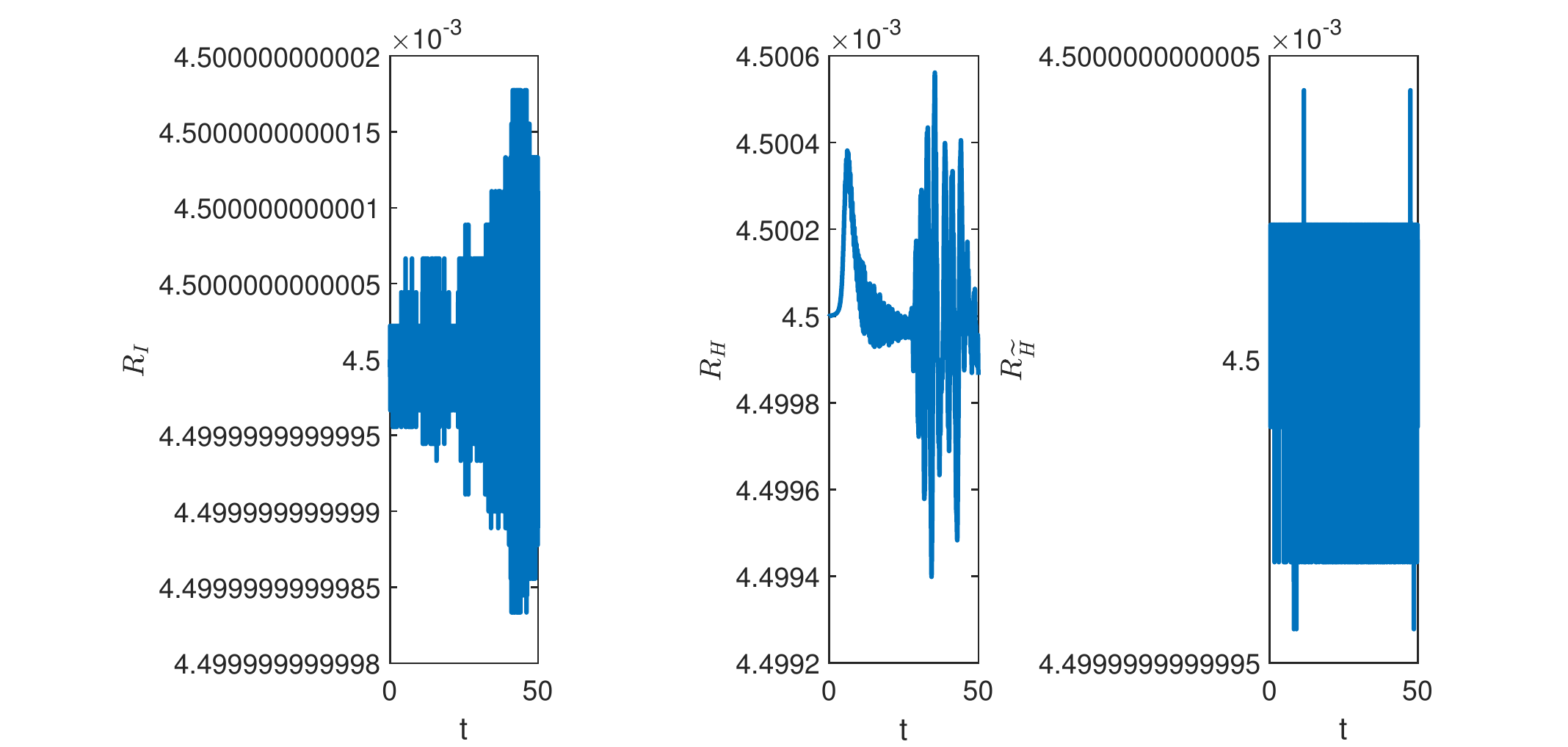}
\caption{Damped Burger's equation with Kahan's method \eqref{Two_step_expkahan}:  residuals of linear conformal invariant (left),   exact Hamiltonian  (middle) and modified Hamiltonian (right).\label{balanceKahan}}
\end{figure}

The solution profiles and  the residual errors are shown, computed with linearly implicit two-step Kahan's method \eqref{Two_step_expkahan}  are shown in Figure \ref{burgerssolkahan} and in Figure \ref{balanceKahan}. The solution deteriorates as time increases, the mass is not preserved as for the EK in Figure~\ref{hamburgers}. These indicate, that the standard linearly implicit methods are not adequate for the damped Hamiltonian Burger's equation.  

\subsection{Damped Korteweg-de-Vries (KdV) equation}

We consider the following  damped KdV equation  \cite{Bhatt21}
\begin{equation}\label{dkdv}
u_t = \alpha (u^2)_x  +  \rho u_x +  \nu u_{xxx} -2\gamma u.
\end{equation} 
In Hamiltonian form, it reads as  
\begin{equation}\label{kdv}
u_t = \frac{\partial}{\partial x} \frac{\delta {\mathcal H}}{\delta u} -2\gamma u, \quad {\mathcal H}(u) = \int \left(\frac{\alpha}{3}u^3 + \frac{\rho}{2}u^2 -\frac{\nu}{2}u_x^2     \right)dx.
\end{equation}
The damped KdV equation was solved with the projected exponential Runge-Kutta methods  \cite{Bhatt21} and   with the conformal multisymplectic method  \cite{Guo19}. It  possesses linear  $I_1 = \int u dx $ and quadratic conformal invariants  $I_2 = \int u^2 dx $  \cite{Bhatt21,Guo19}.

The semi-discrete form of \eqref{kdv} given as   
$$
\dot{{\bm u}} = \alpha D_1 \left( {\bm u} \odot  {\bm u} \right) + \rho D_1 {\bm u} + \nu D_3 {\bm u} - 2 \gamma {\bm u},
$$
where the matrix $D_3=D_1D_2$ approximates the third order derivative $\partial_{xxx}$.  The discrete Hamiltonian has the form
$$
H({\bm u}^n)  = \int \left (  \frac{\alpha}{3}( {\bm u}^n )^3  + \frac{\rho}{2}  ( {\bm u}^n )^2    -   D_1\frac{\nu}{2}( {\bm u}^n )^2  \right ) dx. 
$$

The exponential integrators in Section~\ref{sec2}, yield the following schemes:
\begin{itemize}

\item CIMP: 

\begin{align*}
\frac{e^{X_{1}}{\bm u}^{n+1} - e^{X_0}{\bm u}^n  }{\Delta t }
 = &  \alpha D_1 \frac{e^{X_0}{\bm u}^n \odot e^{X_0}{\bm u}^n + e^{X_{1}}{\bm u}^{n+1}\odot  e^{X_{1}}{\bm u}^{n+1}  }{2}\\
&  + (\rho D_1 + \nu D_3) \frac{e^{X_0}{\bm u}^n + e^{X_{1}}{\bm u}^{n+1}}{2}
\end{align*}

\item EK:

\begin{align*}
\frac{e^{X_{2}}{\bm u}^{n+2} - e^{X_0}{\bm u}^n  }{2\Delta t }  = & \alpha D_1\frac{e^{X_1}{\bm u}^{n+1}\odot  e^{X_0}{\bm u}^n  + e^{X_2} {\bm u}^{n+2}\odot  e^{X_1}{\bm u}^{n+1} }{4  }\\
&  (\rho D_1 + \nu D_3)\frac{e^{X_0}{\bm u}^n + 2e^{X_{1}}{\bm u}^{n+1} + e^{X_{2}}{\bm u}^{n+2}     }{4}.
\end{align*}

\end{itemize}

\begin{figure}[th]
\centering
\includegraphics[width=.55\columnwidth]{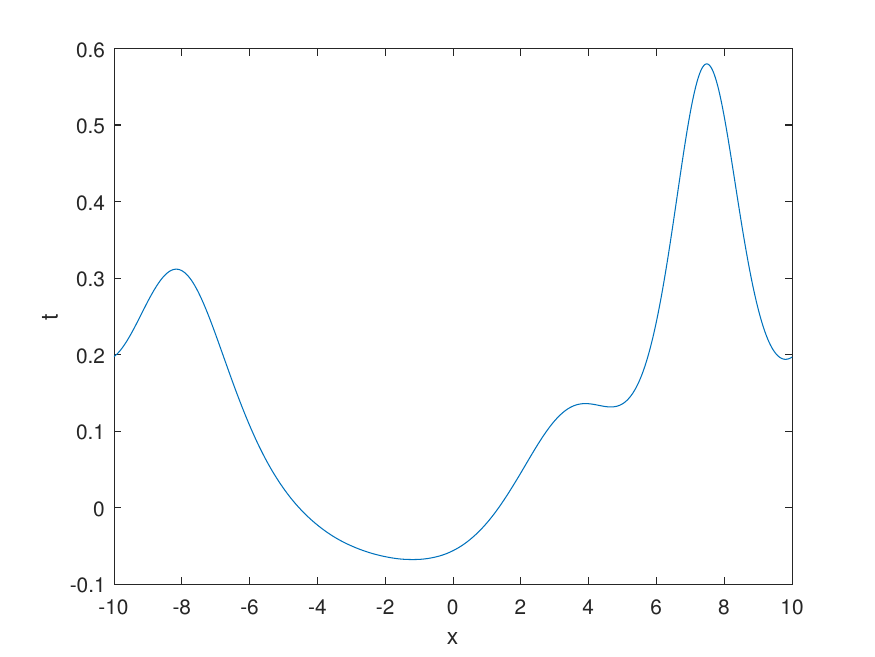}
\caption{Solution of the damped KdV equation at final time.\label{kdvsolution}}
\end{figure}

We consider the KdV equation on  $x\in [-10,10]$   \cite{Bhatt21} with the parameter values
$$
\nu =-10^{-5}, \quad \gamma = 10^{-2}, \quad \alpha = -3/8, \quad \rho = -10.
$$
The mesh sizes are $\Delta x=0.0808$ and $\Delta t=0.009$ with the final time $T=50$. The initial condition is given with the Gaussian wave profile 
$
u(x,0) = \frac{2e^{-2x^2}}{\sqrt{2\pi}}.
$
The numerical
solutions develop almost a vertical front, in form a series of wave-trains in Figure~\ref{kdvsolution} as in \cite{Bhatt21}.

\begin{figure}[th]
\centering
\includegraphics[width=.8\columnwidth]{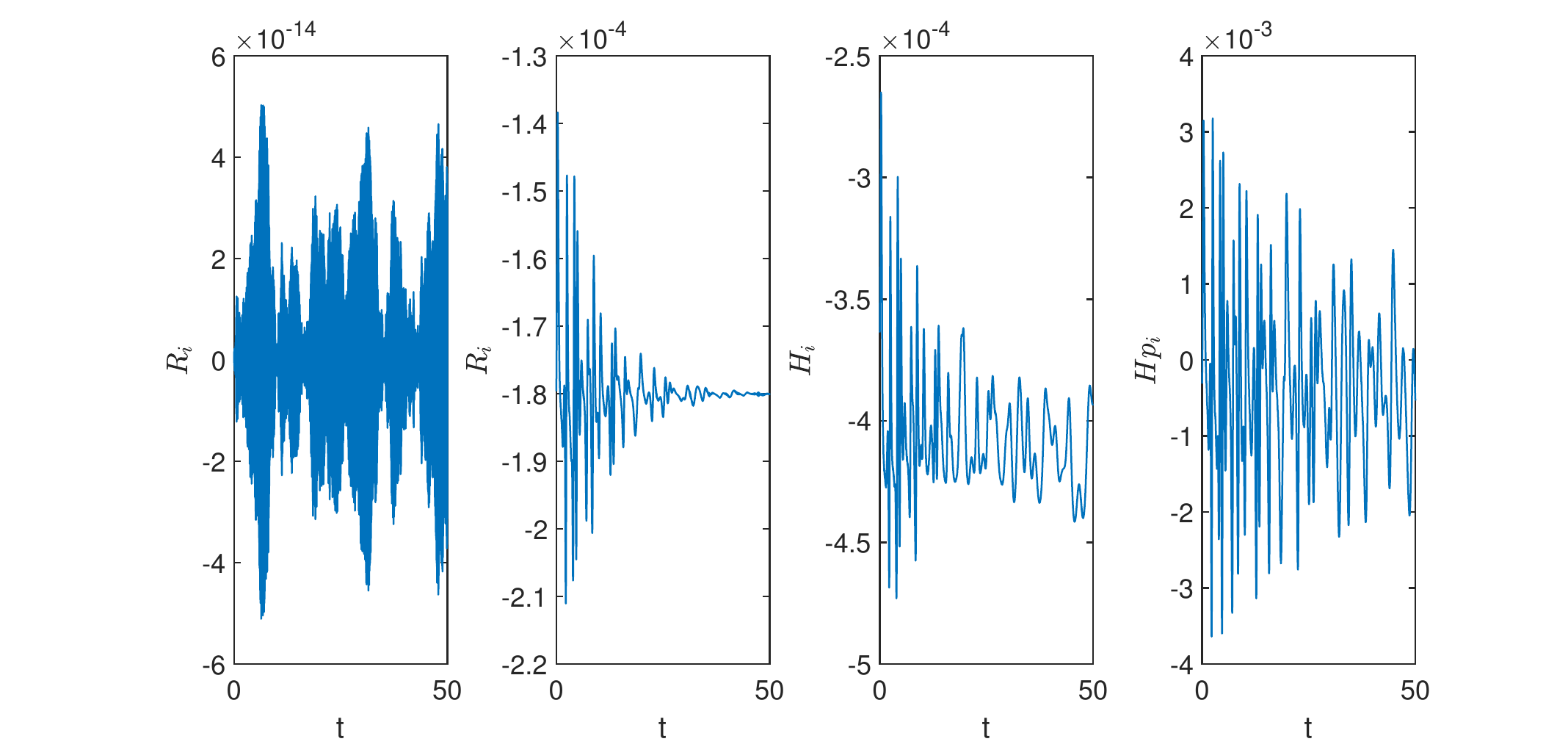}
\caption{Damped KdV equation (EK): residuals of linear and quadratic conformal invariants (left),  Hamiltonian  (middle) and modified Hamiltonian (right).\label{kdvresidualK}}
\end{figure}

The EK preserves the linear conformal invariant up to machine precision,  whereas the error in the residual of the quadratic conformal invariant is diminishing  as time progresses in Figure \ref{kdvresidualK}.

\subsection{Damped nonlinear Schr{\"o}dinger (NLS) equation}

We consider the following damped NLS equation \cite{Cai20,Fu16,Jiang18,Moore13}
\begin{equation} \label{dnls}
	i \psi_t  = -\psi_{xx} -\alpha |\psi |^2\psi - i \frac{\gamma}{2}\psi,
\end{equation}
where $\alpha > 0$ is a constant parameter, and $\gamma >0 $ is a damping coefficient. 
The equation \eqref{dnls} can be written through decomposing $\psi = u + i v$
in real and imaginary components as
\begin{equation} \label{dnls1}
\begin{aligned}
u_t & = -v_{xx} -\alpha (u^2 + v^2) v- \frac{\gamma}{2}u, \\
v_t & =  u_{xx} + \alpha (u^2 + v^2) u- \frac{\gamma}{2}v. \\
\end{aligned}
\end{equation}
The system  \eqref{dnls1}  can be recast into a  damped Hamiltonian system
\begin{equation*}
\begin{pmatrix}
u_t \\ v_t
\end{pmatrix}
=  \begin{pmatrix}
0 &  -1 \\ 1 & 0
\end{pmatrix}
\begin{pmatrix}
\frac{\delta {\mathcal H}}{\delta u }  \\ \frac{\delta {\mathcal H}}{\delta v }
\end{pmatrix}
- \frac{\gamma}{2} \begin{pmatrix}
u \\ v
\end{pmatrix},
\end{equation*}
with the Hamiltonian
$$
\quad {\mathcal H}(u,v)  = \int \left( \frac{\alpha}{4} (u^2 + v^2)^2  -\frac{1}{2} (u_x^2 + v_x^2) \right)dx.
$$

Semidiscretization with finite differences gives the following ODE system 
\begin{equation} \label{dnlsssem}
\begin{aligned}
\dot{{\bm u}} &= - D_2 {\bm v}  -\alpha  \left( {\bm u} \odot  {\bm u} + {\bm v} \odot  {\bm v}  \right)\odot  {\bm v} - \frac{\gamma}{2}{\bm u}, \\
\dot{{\bm v}} &=  D_2 {\bm u}  + \alpha  \left( {\bm u} \odot  {\bm u} + {\bm v} \odot  {\bm v} \right)\odot  {\bm u} - \frac{\gamma}{2}{\bm v},
\end{aligned}
\end{equation}
with the discrete Hamiltonian 
$$
H({\bm u}^n,{\bm v}^n)  = \int \left( \frac{\alpha}{4} ((({\bm u}^n )^2 +  ({\bm v}^n )^2 ))^2  -\frac{1}{2} (D_1({\bm u}^n )^2 +  D_1 ({\bm v}^n )^2 ) \right)dx.
$$

The damped NLS equation \eqref{dnls} has two quadratic conformal invariants;  the mass  $I_1 = \int (u^2 +v^2) dx $ and the momentum $I_2 = \int (v_xu - u_xv) dx $. 
 We solve it with the implicit EAVF and with the linearly implicit LIE integrators.

\begin{itemize}
\item Applying the   EAVF method \eqref{confavf} to  \eqref{dnlsssem} gives the following scheme

\begin{equation*} 
\begin{aligned}
\frac{e^{X_{1}}{\bm u}^{n+1} - e^{X_0}{\bm u}^n  }{\Delta t } & = - \frac{1}{2} D_2 ( e^{X_{1}}{\bm v}^{n+1} + e^{X_0}{\bm v}^n)  - \alpha \int_0^1 \left( {\bm u}_{\xi}^n \odot {\bm u}_{\xi}^n + {\bm v}_{\xi}^n \odot {\bm v}_{\xi}^n \right)\odot {\bm v}_{\xi}^n  d\xi,\\
\frac{e^{X_{1}}{\bm v}^{n+1} - e^{X_0}{\bm v}^n  }{\Delta t } & = \frac{1}{2} D_2 ( e^{X_{1}}{\bm u}^{n+1} + e^{X_0}{\bm u}^n)  + \alpha \int_0^1 \left( {\bm u}_{\xi}^n \odot {\bm u}_{\xi}^n + {\bm v}_{\xi}^n \odot {\bm v}_{\xi}^n \right)\odot {\bm u}_{\xi}^n  d\xi,
\end{aligned}
\end{equation*}
where
$$
{\bm u}_{\xi}^n = \xi e^{X_{1}}{\bm u}^{n+1} + (1-\xi)e^{X_{0}}{\bm u}^{n}, \qquad
{\bm v}_{\xi}^n = \xi e^{X_{1}}{\bm v}^{n+1} + (1-\xi)e^{X_{0}}{\bm v}^{n},
$$
and
$
 \quad X_{1} = \frac{\gamma\Delta t}{4}, \quad X_{0}= -\frac{\gamma\Delta t}{4}.
$
\item The LIE \eqref{polscheme} gives the scheme
\begin{small}
\begin{equation*} 
\begin{aligned}
\frac{e^{X_{2}}{\bm u}^{n+2} - e^{X_0}{\bm u}^n  }{2\Delta t }
& =   - \frac{1}{2} D_2 \left( e^{X_{2}}{\bm v}^{n+2} + e^{X_0}{\bm v}^n  \right) \\
& \qquad - \frac{\alpha}{2}\left(e^{X_1}{\bm u}^{n+1}\odot e^{X_1}{\bm u}^{n+1}  + e^{X_1}{\bm v}^{n+1}\odot e^{X_1}{\bm v}^{n+1} \right)\odot \left( e^{X_{2}}{\bm v}^{n+2} + e^{X_{0}}{\bm v}^{n}\right), \\
\frac{e^{X_{2}}{\bm v}^{n+2} - e^{X_0}{\bm v}^n  }{2\Delta t }
& =   \frac{1}{2} D_2 \left( e^{X_{2}}{\bm u}^{n+2} + e^{X_0}{\bm u}^n \right) \\
& \qquad + \frac{\alpha}{2}\left(e^{X_1}{\bm u}^{n+1}\odot e^{X_1}{\bm u}^{n+1}  + e^{X_1}{\bm v}^{n+1}\odot e^{X_1}{\bm v}^{n+1} \right)\odot \left( e^{X_{2}}{\bm u}^{n+2} + e^{X_{0}}{\bm u}^{n}\right),
\end{aligned}
\end{equation*}
\end{small}
where
$
X_{2} = \frac{\gamma\Delta t}{2},
 \quad X_{1} = 0, \quad X_{0}= -\frac{\gamma\Delta t}{2}.
$
The polarized Hamiltonian has the form 
\begin{small}
\begin{eqnarray*}
\widetilde{H}({\bm u}^n, {\bm u}^{n+1}, {\bm v}^n, {\bm v}^{n+1} )   & =  &\int \left (  \frac{\alpha}{4}(( {\bm u}^n )^2({\bm u}^{n+1})^2  + ( {\bm v}^{n} )^2({\bm u}^{n+1})^2 +  {\bm u}^n {\bm v}^n +
{\bm u}^{n+1} {\bm v}^{n+1}   ) \right .\\
     & & \left .-  D_1\frac{1}{2}(( {\bm u}^n )^2 + ( {\bm u}^{n+1} )^2) - D_1\frac{1}{2}(( {\bm v}^n )^2 + ( {\bm v}^{n+1} )^2)    \right ) dx.
\end{eqnarray*}
\end{small}
\end{itemize}

For the numerical experiment, we consider the damped NLS equation \cite{Cai22,Jiang18} on $x\in[-25,25]$  with the mesh sizes $\Delta x = 50/1024$ and $\Delta t = 0.001$, and the target time is $T=10$. We fix the parameter $\alpha =2$, and take 
the initial condition as
$
\psi(x,0) = \text{sech}(x)e^{2ix}
$.

\begin{figure}[th]
\centering
\includegraphics[width=.50\columnwidth]{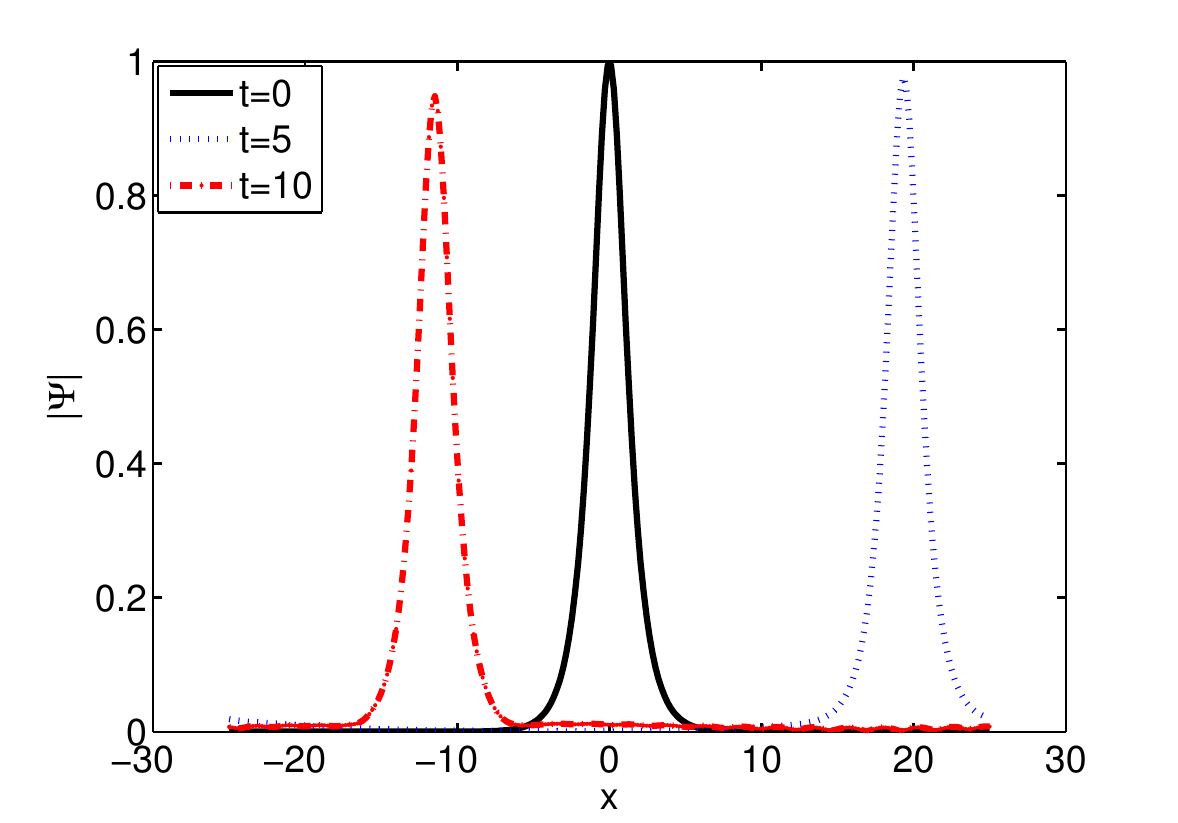}
\caption{Solution of the damped NLS equation.\label{nlssolution}}
\end{figure}

Figure~\ref{nlssolution} shows solutions at initial and final times by using the EAVF and  
LIE schemes with $\gamma=5e-4$.  The damped solitary wave
is traveling from left to right as
required, by  preserving  the phase space
structure.

\begin{figure}[th]
\centering
\includegraphics[width=.45\columnwidth]{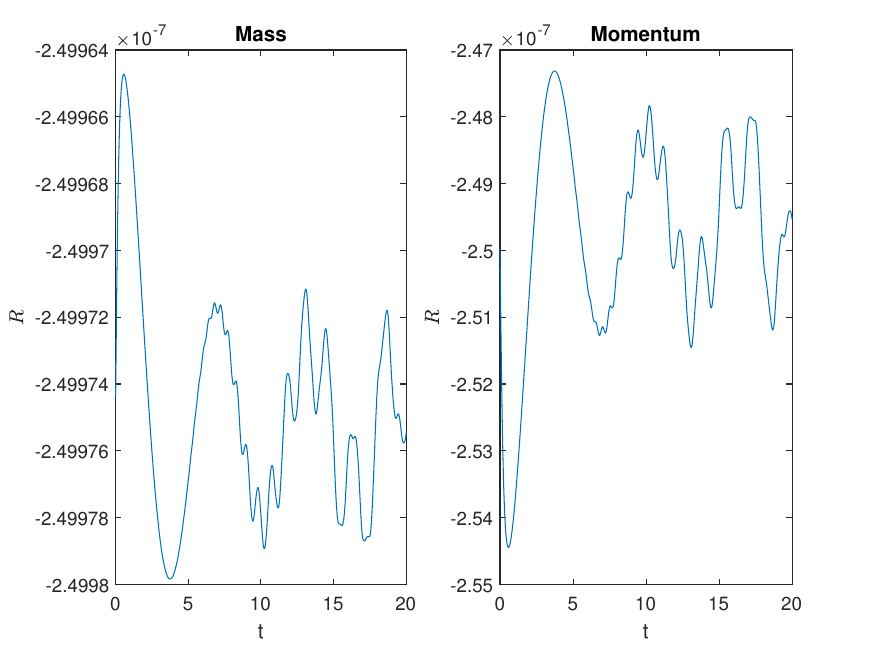}
\includegraphics[width=.45\columnwidth]{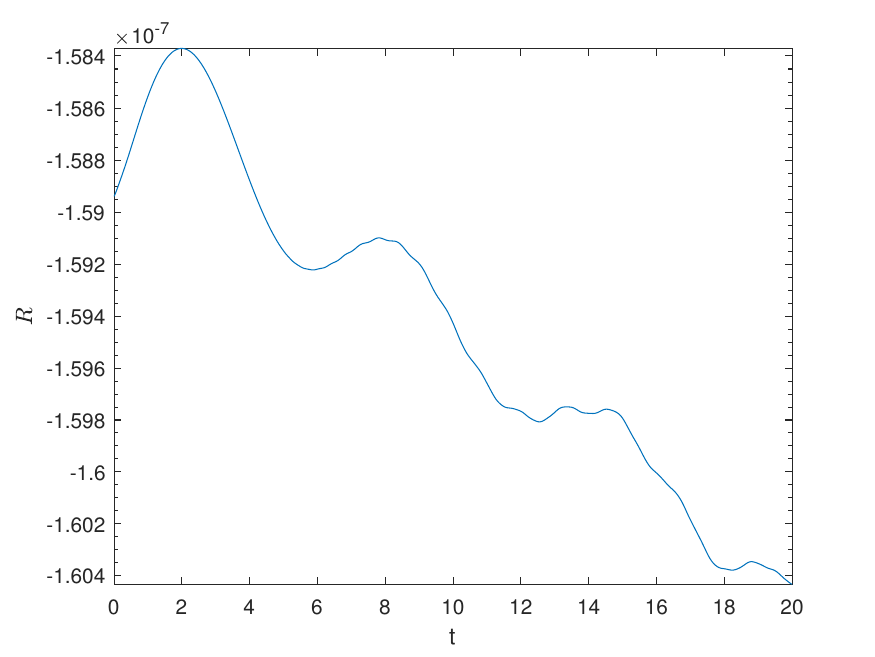}
\caption{Damped NLS equation with $\gamma = 5e-4$ (EAVF): residuals of mass  (left),  momentum  (middle) and  Hamiltonian (right).\label{nlsconformalresidualA}}
\end{figure}

\begin{figure}[H]
\centering
\includegraphics[width=.45\columnwidth]{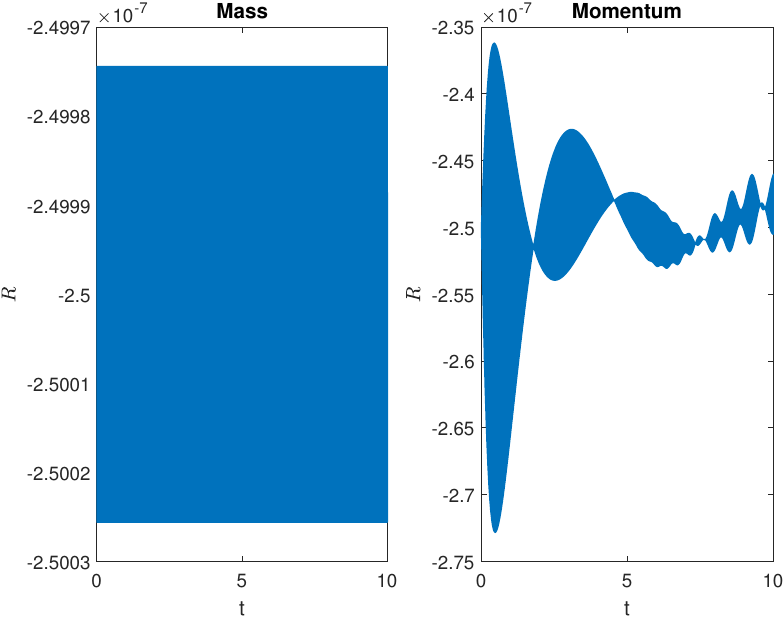}
\includegraphics[width=.45\columnwidth]{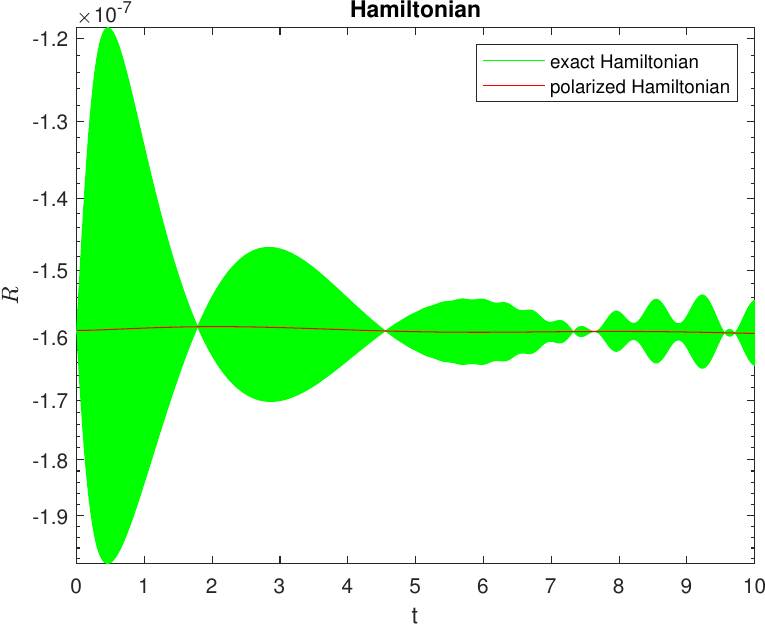}
\caption{Damped NLS equation with $\gamma = 5e-4$ (LIE): residuals of mass (left),  momentum  (middle) and Hamiltonians (right).\label{nlsconformalresidualP2}}
\end{figure}

In Figure \ref{nlsconformalresidualP2}, residuals of the energy balance of the quadratic conformal invariants and the Hamiltonian are plotted. The EAVF and LIE schemes do not preserve the dissipative rate of the invariants exactly, whereas the error in the residuals for LIE is lower than the EAVF. In Figure~\ref{nlsconformalresidualP2}, errors in the residuals are small for small damping factors $\gamma=5e-4$ plotted, whereas the error in residual of the polarized Hamiltonian is exactly preserved.

The CPU time needed for the solution of the systems with LIE is $55.1$ seconds, and with the EAVF $70.4$ seconds, which supports that the linearly implicit exponential integrators are computationally more advantageous than the implicit exponential integrators.

\section{Conclusions}
\label{sec:conc}

Linearly implicit exponential integrators preserve linear conformal invariants exactly and preserve the quadratic invariants,  cubic and quartic Hamiltonians more accurately when the damping coefficient is small. 
Symmetry of the linearly implicit exponential integrators guarantees stable long-time behavior of the solutions. 
Compared with the implicit exponential integrators, linearly implicit exponential integrators  show a lower computational cost as illustrated by the damped NLS equation in long-term integration. 
The computational advantages of the linearly implicit exponential integrators would be more significant for higher dimensional PDEs, which is the subject of future research as well as an extension to time-dependent damping. 


\end{document}